\documentclass[11pt]{article}

\usepackage[hyphens]{url}
\usepackage{mathptmx}      
\usepackage{helvet}        
\usepackage{courier}
\usepackage{graphicx}
\usepackage{amssymb}
\usepackage[margin=1in]{geometry}
\usepackage{subfig}       
\usepackage{type1cm}        
\usepackage{makeidx}         
\usepackage{graphicx} 
\usepackage{float}      
\usepackage{multicol} 
\usepackage{subfig}       
\usepackage{type1cm}
\usepackage{tabularx}        
\usepackage{lipsum}
\usepackage{hyperref}   
\usepackage[justification=centering]{caption}   
\newcommand\blfootnote[1]{%
  \begingroup
  \renewcommand\thefootnote{}\footnote{#1}%
  \addtocounter{footnote}{-1}%
  \endgroup}
\newcommand\T{\rule{0pt}{2.6ex}} % Top strut 
\newcommand\B{\rule[-1.2ex]{0pt}{0pt}} % Bottom strut

\begin{document}

\title{Gender Differences in First Jobs for New US PhDs\\ in the Mathematical Sciences}
\author{Marie A. Vitulli \\
University of Oregon\\
Eugene, OR 97403-1222}
%\institute{\textcopyright  \hspace{.05 cm} 2017 \hspace{.05cm} Marie A. Vitulli  \at University of Oregon,
%Eugene, OR 97403-1222, \email{vitulli@uoregon.edu}}
%\date{today}

\maketitle

\abstract{ 
We take a long term look\blfootnote{\textcopyright  \hspace{.05 cm} 2017 \hspace{.05cm} Marie A. Vitulli} at initial employment trends for new doctorates with an eye towards gender, citizenship, and gender and citizenship differences by analyzing data from 
1991--2015 AMS-ASA-IMS-MAA-SIAM Annual Surveys.  The data show that the unemployment rate for women has been equal to or lower than the
  rate for men during most of the last quarter century.  The one exception is that between 2001 and
   2015 the unemployment rate for women who are not U.S. citizens was higher than the rate for
    non-citizen men. 
  The unemployment rates are higher
 for males who are U.S. citizens than for non-citizen males in the last fifteen years, a puzzling 
 trend. The data show that men from all pure math programs\footnote{When we speak of pure mathematics departments we exclude departments in applied mathematics, statistics, and biostatistics.} are considerably more likely than 
 women to take jobs at the top-ranking and top-producing math departments. The data show
 women take jobs at departments in which the highest degree is a bachelor's degree at 
much higher rates and men take jobs in business and industry at considerably higher 
 rates.   We also find that men from the top-ranking or 
 top-producing doctoral programs tend to be more likely to take jobs at academic 
 institutions or research institutes at least on a par with their degree-granting institutions.
 }
\section{Introduction}

%In the 1990s some mathematicians questioned whether affirmative action efforts were skewing the 
%job market in favor of women.  With this in mind, in 1996, Mary E. Flahive and the current author
% analyzed data from the 1991--95 AMS-IMS-MAA Annual Surveys on initial employment of PhDs in
%  mathematics for possible gender bias 
%in the employment of new  PhD mathematicians \cite{Vitulli Flahive 1997}. In 2010 we analyzed the data for 1996--2008  \cite{Flahive Vitulli 2010}. 
In this study we investigate employment patterns for new PhDs in mathematics between
 1991 and 2015 with an eye toward gender, citizenship and gender $\times$ citizenship\footnote{When social scientists and statisticians study the effects of two independent variables and their interaction they use the mathematical symbol $\times$ to denote testing for an interaction effect.}
 differences  in unemployment rates, patterns of job types, and comparable employment rates, 
 which we  explain in section four. We are concerned about citizenship differences, in large part, because more than half the mathematics PhDs granted by U.S. institutions are awarded to non-U.S. citizens.  The data we analyze comes from the Report on New 
 Doctoral Recipients (from U.S. institutions), which is part of
 the Annual Survey of the Mathematical Sciences published by the American Mathematical Society (AMS).  The Report contains data on jobs both in and outside of academia and in and outside the 
 U.S.  
 More detailed information on the data we used in this study appears at the end of this article.

We raise four questions about employment trends for new PhDs.
\begin{itemize}
\item Is the percentage of women (respectively, U.S. citizen) new PhDs increasing?
\item Are there differences in initial unemployment rates due to gender, citizenship, or 
\mbox{gender $\times$ citizenship}, for new PhDs from U.S. institutions?
\item Are there  differences in the type of employment by gender?
\item  With regard to academic jobs, 
are men and women equally likely be employed by departments whose ranking is at least comparable to the degree-granting department?  
\end{itemize}

We will see that women are still not getting more than their share of the jobs, but the same differences that Flahive and the current author observed in previous studies \cite{Vitulli Flahive 1997} and \cite{Flahive Vitulli 2010} exist today.  For example, women were employed at academic institutions whose highest degree in mathematics is a bachelorÕs degree at a substantially higher rate than men and men were employed in business and industry at a considerably higher rate than women.  

\section{Percentages of Women and U.S. Citizens Among New PhDs}
During the last quarter century women received an average of 28.8\% of the mathematics 
doctorates from U.S. departments.  During this period U.S. citizens received 46.0\% of the 
doctorates;  women received 28.5\% of the doctorates that were awarded to U.S. citizens.  We 
now look at gender $\times$ citizenship differences. We note that the PhDs whose citizenship was
 unknown at the time of the surveys appear in the All PhDs column, but do not appear in the US or
  Non US columns.

\medskip
%\begin{minipage}{\textwidth}
%:  Table of PhDs by Gender and Citizenship
% \begin{center}
 \begin{table}[H]
%\normalsize
\centering
 \caption{Number of PhDs  Granted by Citizenship and Gender 
 (Percentages of all PhDs in \\Two Column Citizenship Group)}
 \label{PhDs by Citizenship and Gender}
%%%  How to explain percentages are for women and men separately

 \begin{tabular}{|c|| c| c ||c|c||c|c|} \hline
& \multicolumn{2}{c||}{US}&\multicolumn{2}{c||}{Non US} & \multicolumn{2}{c|}{All PhDs} \\ \hline \hline
Period& F &M &F &M &F &M \\ \hline
1991--2000 & 1350  & 3622  & 1147  & 4293  & 2597  & 8291   \\
& (27.2\%) & (72.8\%) &  (21.1\%) &  (78.9\%) &  (23.9\%) &  (76.1\%)  \\
2001--2011 & 1950  & 4531 &  2501 & 5187 & 4468  & 9761\\
 & (30.1\%) &  ( 69.9\%) &  (32.5\%) & (67.5\%) & (31.4\%) & (68.6\%)\\
 2012--2015 & 972  & 2548 & 1357  & 2587 & 2330  & 5138 \\
 &  (27.6\%) &  (72.4\%) &  (34.4\%) &  (65.6\%) & (31.2\%) &  (68.8\%)\\
  \hline
 1991--2015 &  4272 & 10701 & 5005 & 12067 & 9395 & 23190 \\
& (28.5\%) & (71.5\%) & (29.3\%) & (70.7\%) & (28.8\%) & (71.2\%) \\
 \hline
 \end{tabular} 
%\caption{PhDs by Gender and Citizenship}
% \end{center}
 \end{table}
%\end{minipage}

\medskip

 Notice that the percentage of new women PhDs was lowest during 1991--2000; during this period the rate for U.S. citizens was considerably higher.  In both 2001--2011 and 2012--2015 non-citizen women received a higher percentage of the degrees than the U.S. citizen women. Among U.S. citizens, the percentage of new women PhDs increased slightly in the second period but fell in the most recent period.   Among the non-citizens, he percentage of new women PhDs  steadily increased.  Non-citizen men received a higher percentage of the degrees than citizens during 1991--2000 but the reverse was true in both later periods.  We find these trends disturbing.
 
\section{Initial Unemployment Rates}

 During the last quarter century, 27,504 out of 32,585 new PhDs  (84.4\%) 
 %checked
 were known to have jobs in or outside of the U.S. by the time the annual survey was conducted in the year in which they received their degrees.  Of the 5,081 new PhDs that didn't report having jobs, 3,464 either remained in the U.S. and had unknown employment status or left the U.S. and didn't report having jobs.\footnote{The AMS divides new PhDs who remain in the U.S. after receiving their degrees without reported employment into three groups: Still Seeking US, Not Seeking US, and Unknown US.  In contrast, all people who 
leave the U.S. after their degrees are reported as 
Unknown Non US.}   
 
We now will take a closer look at initial unemployment rates, focusing on those who remained in the U.S. after receiving their degrees 
and were still seeking employment at the time of the survey.  Following current AMS conventions on unemployment rate calculations, individuals employed outside the U.S. as well as those whose employment status was unknown have been removed from the denominator in the calculation of the unemployment rate.  We also adopt the AMS convention of removing those individuals reported as not seeking employment from the denominator. %We did not consider those new PhDs who were not seeking employment at time of the survey.  
For the entire 25 year period under investigation, the group Not Seeking US accounted for 1.1\% of all new PhDs (1.6\% of the females and 0.9\% of the males).  As pointed out in the 2015 survey, these conventions increase the unemployment rate from the rates reported prior to the these adjustments.

During the 1991--2015 time period, 1,598  of 24,841 (6.4\%) of the these new PhDs were still seeking employment at the time of the annual survey; the rate was 0.6\% higher for non-U.S. citizens than for citizens.
  
  Table \ref{PhDs Still Seeking by Citizenship}  breaks down the unemployment rate by citizenship over three time periods in the
   last quarter century. In all tables in this section the \textit{All PhDs} column includes new PhDs 
   whose citizenship was unknown at the time of the survey.

 %%% Both Tables Checked for Accuracy %%%
%: Table of PhDs Still Seeking Employment by Citizenship
\medskip
%\begin{minipage}{\textwidth}
% \begin{center}

\begin{table}[h]
%\normalsize
\centering
\label{PhDs Still Seeking by Citizenship}
\caption{Number of New PhDs in the US Still Seeking Employment
 by Citizenship\\ (Percentages of all New PhDs in Column Cohort)}

\begin{tabular}{|c|| c| c ||c|} \hline
Period & US & Non US   & All PhDs \\ \hline
1991--2000 & 332 (7.6\%) & 388 (11.2\%) & 749 (9.3\%)  \\
2001--2011 & 246 (4.5\%) & 240 ( 4.5\%) & 486 ( 4.5\%) \\
 2012--2015 & 217 (7.1\%) & 146 (5.3\%) & 363 (6.2\%) \\
 \hline
 1991--2015 &795(6.1\%) & 774 (6.7\%) & 1598 (6.4\%)\\
 \hline
 \end{tabular} 
%\caption{PhDs Still Seeking Employment}

\label{PhDs Still Seeking by Citizenship}
% \end{center}
% \end{minipage}
\end{table}
 \medskip
% \normalsize
Notice that during 1991--2000 
 the unemployment rate for non-U.S. citizens was 
% 3.6 percentage points  or 
 47.4\% higher than the rate for U.S. citizens.  The unemployment rate for 2001--2011 was 
 the same for both U.S. citizens and non-U.S. citizens and was lower than the rate for the preceding decade.   The overall unemployment rate increased in 2012--2015. Looking
at citizenship differences, the rate for  for non-U.S. citizens was 
%1.8 percentage points or 
25.5\% lower than the rate for U.S. citizens.  The early disadvantage for 
  non-U.S. citizens was reversed by the end of the study.    It is 
  disappointing that the unemployment rate is 
  inching back up in recent years after a decline during 2000--2011.  The
  reader is reminded that due to the change in groupings, the cycles of new PhDs are unequal in length and 
  hence  these percentage differences are suggestive, but not directly comparable.  Over the entire 
  1991--2015 period, the unemployment rate for non-U.S. citizens was 0.6 percentage points higher, which represents  8.6\%   of the unemployment rate for citizens.
  
 In Table \ref{PhDsStillSeekingbyCitizenshipandGender}, we look at three separate time periods and break down the unemployment rates by gender and citizenship. For U.S. citizens, the rate of those still seeking employment is  lower for women during each time period whereas for non-U.S. citizens the rate is at least at high for women during each time period.  

\medskip
%\begin{minipage}{\textwidth}
%:  Table of PhDs Still Seeking Employment by Citizenship and Gender
 \begin{table}[H]
 \centering
\caption{Number of New PhDs Still Seeking Employment
 by Citizenship and Gender \\(Percentages of all New PhDs in Column Cohort)}
 \label{PhDsStillSeekingbyCitizenshipandGender}
%%%  How to explain percentages are for women and men separately
\small
 \begin{tabular}{|c|| c| c ||c|c||c|c|} \hline
& \multicolumn{2}{c||}{US}&\multicolumn{2}{c||}{Non US} & \multicolumn{2}{c|}{All PhDs} \\ \hline \hline
Period& F &M &F &M &F &M \\ \hline
1991--2000 & 61 (5.1\%) & 271 (8.5\%) & 84 (11.2\%) & 304 (11.2\%) &154 (7.7\%) & 595 (9.8\%)   \\
2001--2011 & 50 (2.9\%) & 196 ( 5.1\%) &  89 (4.8\%) & 151 (4.3\%) & 139 (3.9\%) & 347 (4.7\%)\\
 2012--2015 & 41 (4.7\%) & 176 (8.0\%)& 61 (6.1\%) & 85 (4.8\%) & 102 (5.5\%) & 261 (6.6\%)\\
 \hline
 1991--2015 &  152 (4.0\%) & 643 (7.0\%) & 234 (6.5\%) & 540 (6.7\%) & 395 (5.3\%) & 1203 (6.9\%) \\
 \hline
 \end{tabular} 
%\caption{PhDs Still Seeking Employment}

% \end{center}
%\end{minipage}
\end{table}
\normalsize
\medskip

Looking at gender $\times$ citizenship differences, the disadvantage for non-U.S. citizens from 
1991--2000 is more pronounced for women;  the rate for female non-U.S. citizens  was more than
 double the rate for citizens. The disadvantage for non-U.S. citizens lessened during 2001--2011.   
 Notice that for male PhDs during the 2012--2015 period, the unemployment rate for citizens was 3.2 percentage points higher than for non-U.S. citizens, which is
67.3\% of the unemployment rate for non-U.S. citizens.

%: Types of Employment
 \section{Types of Employment}
 
 In this section, we look at the types of employment obtained by the various groups of new PhDs 
We start with general observations with an eye towards gender differences and then break down our observations by public/private degree-granting institutions and citizenship.  As we shall see, some of the trends we saw in 
the earlier studies persist: women take jobs at institutions in which the highest degree is a 
bachelor's degree at substantially higher rates and 
 men take jobs in business and industry at somewhat higher rates. The AMS changed the annual survey reporting groupings in 2012 so we present data separately for 2012 -- 2015. Beginning in 2012 the top-ranking Group I departments were replaced by the top-producing Public Large and Private Large departments. The All Others row in all tables in this section include new PhDs who accepted jobs in statistics, biostatistics or applied math departments, outside the U.S., those who were still seeking or not seeking employment, as well as those whose employment status was unknown at the time of the survey. 

% Although our current analysis is similar to the earlier one, there are some differences.  In the current study the data on new PhDs from Group Va departments are included.  Also, over the years there have been some changes in the recording of data, principally the 1996 change in the groupings of doctorate-granting institutions already discussed. Finally, in Table 2 we report comparable employment percentages for all new PhDs rather than for U.S. citizens only.
\subsection{Pure Mathematics Doctorates:All Departments}

   Let's first group together all the new PhDs who received degrees from pure mathematics 
   programs, that is, Ph.Ds. from Groups I -- III departments during 1991--2011 and from
   Public Large/ Medium/Small and Private Large/Small institutions during 2012--2015.  
   In all, there were 17,753 new PhDs in that cohort, 4,333 (24.4\%) of whom were women and 13,420 (75.6\%) of whom were men.    
  
The first two tables summarize the findings for new PhDs from pure mathematics doctoral 
programs who were employed in the U.S. In Table \ref{Groups_I_III} we report for the entire 1991--2011 time period as well give breakdowns for into 1991--2000 and 2001--2011 periods.  We look at the 2012--2105 cohort separately in Table~\ref{PureMathPhDs}. Since the number of new pure math PhDs who were employed by statistics, biostatistics, applied mathematics, or operations research (Groups IV and V) departments is very small, we do not include separate data rows on employment in those areas in the table; however they are included in the All Others rows as well as the column total and percentage of column cohort calculations.  We remind the reader that in AMS data reports, ``other academic" stands for US academic departments other than pure and applied mathematics departments, biostatistics departments, departments whose highest degree is Bachelor's or Master's degrees, and 2 year colleges. 
%The denominators in the percentage of column cohort calculations also include new PhDs who accepted other nonacademic jobs, those who were still seeking or not seeking employment, as well as those whose employment status was unknown at the time of the survey. This is true for all tables in this section.

\medskip

%\begin{minipage}{\textwidth}
%\begin{center}
\begin{table}[t]
\centering
\caption{Observed Frequencies of First Jobs 
(Percentages of Column Cohort) \\for  Pure Math PhDs 1991--2011}
\label{Groups_I_III}
\smallskip

%: Observed Frequencies of First Jobs for Groups I - III
\footnotesize
%\noindent
\begin{tabular*}{1.05\linewidth}{@{\extracolsep{\fill}} |r ||c| c| c|| c| c| c|| c| c|c|} 
%\begin{tabular*}{1.0\linewidth}  |r |c| c| c| c| c| c| c| c|c|} 

 %\begin{tabular}{|r||c|c|c||c|c|c||c|c|c|} 
 \hline
& \multicolumn{3}{c||}{1991--2000}&\multicolumn{3}{c||}{2001--2011} &\multicolumn{3}{c |}{1991--2011}\\ \hline \hline
  Emp Type&F&M&All&F&M&All &F&M&All  \\ \hline
Gr I& 199  & 849  &1048  & 335  & 1216  & 1551 & 534 & 2065 & 2599  \\
& (11.2\%) & (13.3\%) & (12.9\%) & (13.1\%) & (17.3\%) & (16.2\%) &(12.3\%) & (15.4\%) & (14.6\%) \\  \hline
Gr II & 88  & 326  & 414  & 193 & 559 & 752 & 281 &885 &1166  \\
& (4.9\%) & (5.1\%) & (5.1\%) &  (7.6\%) &  (7.9\%) &  (7.8\%) & (6.5\%) & (6.6\%) & (6.6\%) \\ \hline
Gr III & 99  & 251 & 350  &120  & 263  & 383 &219 & 514 & 733  \\
&  (5.6\%) & (3.9\%) &  (4.3\%) &(4.7\%) &  (3.7\%) & (4.0\%) & (5.1\%) & (3.8\%) & (4.1\%) 
 \\ \hline
Masters &161 & 416  &577  & 202  &351  &553  &363 & 767 & 1130 \\
&(9.0\%) & (6.5\%) &(7.1\%) &  (7.9\%) &(5.0\%) & (5.8\%) &(8.4\%) & (5.7\%) & (6.4\%)  \\ \hline
Bachelors &370   & 787   & 1157   &  483 &877  & 1360 & 853 & 1664 & 2517 \\
& (20.8\%)  & (12.3\%)  & (14.2\%)  &   (18.9\%)& (12.4\%)  & (14.2\%) & (19.7\%) & (12.4\%) & (14.2\%) \\ \hline
2Yr & 38  &110  &148   & 61  &150  & 211 & 99 & 260 & 359 \\
& (2.1\%) & (1.7\%)  &(1.8\%)  & (2.4\%)  & (2.1\%)  &  (2.2\%) & (2.3\%) & (1.9\%) & (2.0\%)  \\ \hline
Oth Acad &  61 & 239   & 300   &190  &407  & 597 &251 & 646 & 897  \\
 &   (3.4\%)&  (3.8\%)  & (3.7\%)  & (7.4\%)  & (5.8\%)  &  (6.2\%)  & (5.8\%) & (4.8\%) & (5.1\%)\\ \hline
Res Inst &  32 & 176  & 208 & 70 & 134 & 204  & 102 & 310 & 412 \\
& (1.8\%)& (2.8\%) & (2.6\%) &  (2.7\%) &  (1.9\%) &  (2.1\%) & (2.4\%) & (2.3\%) & (2.3\%)\\ \hline
Govt & 36  &141 & 177  & 90 & 202  & 292 &126 & 343 & 469 \\ 
&  (2.0\%)  & (2.2\%) &  (2.2\%) & (3.5\%)&  (2.9\%) & (3.0\%) & (2.9\%) & (2.6\%) & (2.6\%)\\ \hline
Bus/Ind & 143  & 754 &897 & 205 &692 &897 &348 & 1446 & 1794  \\ 
&  (8.0\%) &  (11.8\%) &(11.0\%) &  (8.0\%)& (9.8\%) & (9.3\%) & (8.0\%) & (10.8\%) & (10.1\%) \\ \hline
All Others & 555 & 2324 & 2879 &602 & 2196 & 2798 & 1157 & 4520 & 5677 \\
& (31.1\%) & (36.5\%) & (35.3\%) &(23.6\%) & (31.2\%) & (29.2\%) & (26.7\%) & (33.7\%) & (32.0\%) \\
\hline
Grand Total &1782  & 6373 & 8155 &2551 & 7047& 9598& 4333 & 13420 & 17753 \\
& (100.0\%) &  (100.0\%) & (100.0\%) & (100.0\%) & (100.0\%) & (100.0\%) & (100.0\%) & (100.0\%) & (100.0\%) \\
\hline
\end{tabular*}
%\label{Groups I-III}
%\end{minipage}
\end{table}
\medskip
\normalsize

%In Table 3b. below we present the data for pure math PhDs during the 2012--2015 period.

We will make some observations by looking at large differences in percentages of men and women who received pure math PhDs and who were employed in various sectors.

\subsubsection{Group I and Public/Private Large Hires} 
Looking at the period 1991 -- 2000 in Table \ref{Groups_I_III} we see a 2.1 point difference in the percentages of men and women who were employed by the top-ranking Group I departments; this 2.1 point difference represents 18.8\% of the percentage of women who were employed by Group I departments. Between 2001 and 2011 the difference was more pronounced. We see a 4.2 point difference in the percentages of men and women who were employed by the top-ranking Group I departments; this 4.2 point difference represents 32.1\% of the percentage of women who were employed by Group I departments.
To put this in perspective, during the same period, 46.9\% of all women pure math PhDs received degrees from Group I institutions and 56.1\% of the men pure math PhDs received degrees from Group I institutions.   Looking only at Group I PhDs, we saw that 24.4\% were women and 75.6\% were men.
 
 \medskip

We now turn our attention to the 2012 -- 2015 cohort of new pure math PhDs.  We will first present the data and then make some observations.

%\begin{minipage}{\textwidth}
%: Observed Frequencies of First Jobs for Public/Private Large, Med, Small
%\begin{center}
%Table 3b.
\begin{table}[h]
\centering
\caption{Observed Frequencies of First Jobs 
(Percentages of Column Cohort) \\for  Pure Math PhDs 2012--2015}
\label{PureMathPhDs}
\smallskip
%\small
 \begin{tabular}{|r||c|c|c|} \hline 
 Employer Type & Female &Male & All \\ \hline \hline
 Public Large &  94   & 332  & 426   \\ 
  &    (7.4\%) &   (9.3 \%) &   (8.8\%) \\ \hline
 Public Medium & 60   & 179   & 239   \\ 
 &  (4.7\%) &  (5.0\%) &  (5.0\%) \\ \hline
 Public Small & 51   & 129  & 180   \\
 & (4.0\%) &   (3.6\%) &   (3.7\%) \\  \hline
 Private Large & 55  & 230   & 285   \\ 
 &  (4.3\%) &   (6.5\%) &   (5.9\%) \\  \hline
 Private Small & 29   & 62   & 91   \\
 &   (2.3\%) &   (1.7\%) &   (1.9\%) \\ \hline
 Masters & 53   & 121   & 174   \\
 &   (4.2\%) &   (3.4\%) &   (3.6\%) \\  \hline
Bachelors & 222   & 354  & 576  \\ 
&  (17.5\%) &   (10.0\%) &   (12.0\%) \\ \hline
2 Year & 35   & 80  & 115   \\ 
&   (2.8\%) &  (2.3\%) &  (2.4\%) \\ \hline
Other Academic & 91   &197   &  288 \\
&   (7.2\%) & (5.5\%)  &    (6.0\%)\\  \hline
 Research Inst & 25  & 87   & 112    \\ 
 &   (2.0\%) &  (2.4\%) &   (2.3\%)  \\  \hline
 Government & 52  & 109   & 161   \\ 
 &   (4.1\%) &  (3.1\%)  &  (3.3\%)  \\ \hline
 Business/Industry & 171  & 577   & 748   \\
 &   (13.5\%) &   (16.2\%) &   (15.5\%)  \\  \hline
 All Others & 328 & 1095 & 1423 \\
 & (25.9\%) & (30.8\%) & (29.5\%) \\ \hline
 Grand Total & 1266 & 3552 & 4818 \\
 & (100\%)& (100\%)& (100\%) \\ \hline
\end{tabular}
%\end{center}

%\end{minipage}
\end{table}
\medskip
\normalsize
 
    Table \ref{PureMathPhDs} shows that between 2012 and 2015 the percentages of men who 
were employed in Public 
    and Private Large departments were noticeably higher than the percentages for women, with the biggest difference occurring in Private Large hires.    At the Public Large departments, 
    we see a 1.9 point difference in the percentages of men and women who were employed by the top-ranking Group I departments; this 1.9 point difference represents 25.7\% of the percentage of 
    women who were employed by the  Group I departments.  At Private Large departments we see a 2.2
     point difference in the percentages of men and women who were employed by the top-ranking Group I departments; this 2.2 point difference represents 51.2\% of the percentage of women who were employed by Public and Private large departments. The latter was the greatest difference we observed.
    
    %  It is tempting to equate the  categories of Public Large and Private Large with the old categories Group I Public and Group I Private, since the top ranked schools tend to produce the most PhDs.

\subsubsection{Group II and Public Medium Hires} There were small differences in the percentages of men and women who were employed by Group II or Public Medium 
departments. In all time periods of this study, the percentages of male pure math PhDs who were employed by Group II or Public Medium departments were slightly higher that the percentages for females. 
 
\subsubsection{Group III Hires} During 1991--2011 we see a 1.3 point difference in the percentages of men and women who were employed by the Group III departments with women having the higher percentage; this 1.3 point difference represents 34.2\% of the percentage of men who were employed by Group III departments.   During 2012 and 2015 the percentage of women who were employed by Public Small departments was  0.4 points higher than the percentage for men; this 0.4 point difference represents 11.1\% of the percentage of men who were employed by Public Small departments.  Over the same period, the percentage of men who were employed by Private Small departments was 0.6 points higher than the percentage for women; this 0.6 point difference represents 35.3\% of the men who were employed by Private Small departments.

\subsubsection{Bachelor's Degree Only Hires} There were striking gender differences in the 
percentages of male and female new pure math PhDs who took jobs in departments in which the 
highest mathematics degree is a  bachelor's degree during all periods of this study. The reader should first refer 
to Table \ref{Groups_I_III} to follow our observations.  Between 1991 and 2000 the percentage of
 women who took 
jobs at  bachelor's-only departments was 8.5 points higher than for men; this 8.5 point difference 
represents 69.7\% of the percentage of men who were employed by bachelor's-only departments. 
Between 2001 and 2011 the difference was slightly less but still large; 
the percentage of women who took these jobs was 6.5 points higher than the percentage of men; this 6.5 point difference represents 52.4\% of the percentage of men who were employed by bachelor's-only departments. Looking at Table \ref{PureMathPhDs}, we see the percentage difference increased during 2012--2015; the percentage of women was 7.5 points higher  than the percentage of men who were employed by bachelor's-only departments; this 7.5 point difference represents 75.0\% of the percentage of men who were employed by bachelor's-only departments.  This is the largest difference we observed for bachelor's-only hires. This gender difference was substantial in every period under analysis.

%:Left off here 4-21-17 
 
\subsubsection{Business and Industry Hires}
  There were considerable gender differences in the percentages of new pure math PhDs who took  jobs in business and industry.  During 1991 and 
 2000, the percentage of men taking jobs in business and industry was 3.8 points
higher than the percentage for women; this 3.8 point difference represents 47.5\% of the percentage of women who took jobs in business and industry.  For
  2001--2011 the percentage of men taking jobs in business and industry was 1.8 points or
22.5\% higher than the percentage for women. During both of these periods, the percentage of women who took jobs in business and industry remained constant at 8.0\%, whereas the percentage for men who  took jobs in business and industry dropped from 11.8\% to 9.8\%; this drop caused the gender difference to lessen.   In the most recent period, 2012--2105, the percentage of men taking jobs in business and industry was 2.7 points or 
20.0\% higher than the percentage for women.  Notice that the percentage of women who took jobs in business and industry between 2012 and 2015 was 13.5\% compared to 8.0\% for both 1991--2000 and 2000--2011. Thus the percentage gender differences observed for jobs accepted in business and industry has diminished in recent years.

%As in earlier years, men were 25.6\% more likely than women to be hired by Public Large and 51.2\% more likely than women to be hired by Private Large departments.  
%Between 2012 and 2015 the percentages of new pure math PhDs taking jobs at bachelor's only departments were 17.5\% for women and 10.0\% for men; the percentage difference during this period is 75.0\%, which is larger than the previous differences. In this period, the percentages of new pure math PhDs hired in business and industry were 13.5\% for women and 16.2\% for men;  men obtained 20.0\% more of the jobs in business  and industry during these recent years.

\subsection{Pure Mathematics Doctorates: Group I Public/Private \\and Public/Private Large Departments}
The top-ranking or top-producing departments award more of the doctorates than any other group of departments.  Between 1991 and 2011, 38.0\% of all new Ph.Ds and 53.8\% of the pure math PhDs came from Group I 
Public or Private departments. Between 2012 and 2015, 33.3\% of all new PhDs and 51.4\% of the pure math PhDs came from 
Public or Private Large departments.  We now look at how they fared in the U.S. job market.  Since Group I wasn't subdivided into Group I Public and Group I Private until 1996 we start Table \ref{GroupIPhDs} in 1996.
%Throwing in the Public Medium institutions, we can account for 
%49.3\% of the new 2012--2015 PhDs We will take a closer look at employment types for PhDs from these 
%institutions.

We see in Tables ~\ref{GroupIPhDs}  and \ref{PublicPrivateLargePhDs} that the gender differences we saw for jobs obtained at Group I departments by all pure math PhDs are greatly reduced when looking only at new PhDs from the top-ranking or top-producing departments. As in Tables \ref{Groups_I_III} and \ref{PureMathPhDs}, the denominators in the percentage of column cohort calculations include new PhDs who accepted other nonacademic jobs, those who were still seeking or not seeking employment, as well as those whose employment status was unknown at the time of the survey.
We continue
 to see very large gender differences in jobs taken at bachelor's-only departments and substantial 
 differences in jobs taken in business and industry, particularly for Group I Public PhDs. Our data also show
  that more than half of the jobs in business and industry between 1991 and 2011 went to Group I 
  PhDs.  We remind the reader that since the number of these PhDs who obtain jobs in 
  statistics, biostatistics, applied mathematics, or operations research is very small we do not
   include data on employment in those areas as rows in our table; however they are included in column total and percentage of column cohort calculations.

\medskip

%\begin{minipage}{\textwidth}
%\begin{center}
\begin{table}[H]
\centering
\caption{Observed Frequencies of First Jobs 
(Percentages of Column Cohort) \\ for  
Group I PhDs 1996--2011}
\label{GroupIPhDs}
\medskip

%: Observed Frequencies of First Jobs for Group I only
%\small
 \begin{tabular}{|r||c|c|c||c|c|c|} \hline
& \multicolumn{3}{c||}{Group I Public}&\multicolumn{3}{c|}{Group I Private} \\ \hline \hline
  Employer Type&Female&Male&All&Female&Male&All  \\ \hline
Group I Public & 143  & 528   &671  & 68    & 251 & 319   \\ 
&   (13.5\%) &  (14.7\%) & (14.4\%) &  (11.9\%)  &  (11.6\%) &  (11.6\%) \\ \hline
Group I Private & 67   & 288  &355  & 102   & 400  & 502  \\ 
&  (6.3\%) &   (8.0\%) &  (7.6\%) &   (17.9\%)  &  (18.5\%) &  (18.3\%) \\ \hline
Group II & 78   & 299   & 377   & 41   & 120  & 161  \\ 
&   (7.4\%) &  (8.3\%) &  (8.1\%) &   (7.2\%) &   (5.5\%) & (5.9\%) \\ \hline
Group III & 37   & 70 & 107   & 8   & 28  & 36  \\
&  (3.5\%) &   (1.9\%) &  (2.3\%) &   (1.4\%) &   (1.3\%) &   (1.3\%) \\ \hline
Masters & 65  &  123  &188   & 11   &37  &48   \\ 
&   (6.1\%) &   (3.4\%) &  (4.0\%) &  (1.9\%) & (1.7\%) &  (1.8\%)  \\ \hline
Bachelors &168    & 289   & 457   & 54  & 118   & 172  \\ 
  & (15.9\%)  &   (8.0\%)  & (9.8\%)  &  (9.5\%)&  (5.4\%)  &   (6.3\%) \\ \hline

2 Year & 13   &47    &60    & 2   &8   & 10  \\
&  (1.2\%) & (4.4\%)  &  (1.3\%)  &  (0.4\%)  &  (0.4\%)  &  (0.4\%) \\  \hline
Other Academic &  47  & 140   &  187    &28   &107    & 135    \\ 
&   (3.4\%)&  (3.9\%)  &    (4.0\%)  &  (4.9\%)  &  (4.9\%)  &  (4.9\%)  \\  \hline
Research Inst &  27  & 72   & 99  & 26  & 86   & 112   \\
 &   (2.5\%)&   (2.0\%) &  (2.1\%) &  (4.6\%) &  (4.0\%) &   (4.1\%) \\  \hline
Government & 38  &96  & 134   & 16 & 40   & 56  \\ 
&   (3.6\%)  & (2.7\%) &   (2.9\%) &   (2.8\%)&   (1.8\%) &   (2.0\%) \\  \hline
Business/Industry & 90  & 396   &486  & 58 & 246  &304    \\ 
&   (8.5\%) &  (11.0\%) &  (10.4\%) &   (10.2\%)&   (11.3\%) &  (11.1\%)  \\  \hline
All Others & 286 & 1246 & 1532 & 157 & 727 & 884 \\
& (27.0\%) & (34.7\%) & (32.9\%) & (27.5\%) & (33.5\%) & (32.3\%) \\ \hline
Grand Total & 1059 & 3594 & 4653 & 571 & 2168 & 2739 \\
& (100.0\%) & (100.0\%) & (100.0\%) & (100.0\%) & (100.0\%) & (100.0\%)  \\ \hline
\end{tabular}
%\end{center}
%\end{minipage}
\end{table}

\normalsize

Notice that for both Group I Public and Group I Private new doctorates, the percentage of women who were employed by bachelor's-only departments is nearly double the percentage for men. 

%: Observed Frequences for Public/Private Large PhDs

%\begin{minipage}{\textwidth}
%\begin{center}
%Table 3d. 
\begin{table}[ht]
\centering
\caption{Observed Frequencies of First Jobs 
(Percentages of Column Cohort)\\ for  
Public/Private Large  PhDs 2012--2015}
\label{PublicPrivateLargePhDs}

 \begin{tabular}{|r||c|c|c||c|c|c|} \hline
& \multicolumn{3}{c||}{Public Large} & \multicolumn{3}{c|}{Private Large} \\ \hline \hline
  Employer Type&Female&Male&Totals&Female&Male&Totals  \\ \hline
  Public Large &  53  & 177   & 230  & 22  &88    &110 \\
  &   (14.3\%) &   (14.2 \%) &  (14.2\%)&   (12.0\%) &  (13.1\%)  &  (12.8\%)\\  \hline
Public Medium & 18   & 47   & 65  &6   & 30   &  36  \\ 
&  (4.9\%) &   (3.8\%) &  (4.0\%) &  (3.3\%) &   (4.5 \%) &  (4.2\%)\\ \hline
 Public Small & 5  & 30   & 35 & 1  &5    &6   \\ 
 &  (1.4\%) &  (2.4\%) &  (2.2\%) &  ( 0.5\%) &   (0.7 \%) &   ( 0.7\%) \\  \hline
Private Large & 17   & 98   & 115   & 31    &103    &134   \\
&  (4.6\%) &  (7.8\%) &  (7.1\%) &  (16.9\%)  &   (15.3\%)  &   (15.6\%) \\   \hline
Private Small & 7 & 19   & 26   & 5   &9    &  14  \\ 
&  (1.9\%) &   (1.5\%) &   (1.6\%) &  (2.7\%)  &   (1.3\%)  &   (1.6\%) \\  \hline
 Masters & 7  & 20   & 27   & 4   & 8   & 12  \\ 
 &  (1.9\%) &  (1.6\%) &   (1.7\%) &   (2.2\%) &   (1.2\%) &   (1.4\%)\\ \hline
Bachelors & 55   & 83   & 138   & 12  & 22   & 34   \\ 
&  (14.9\%) &   (6.6\%) &   (8.5\%) &   (6.6\%) &  (3.3\%) &   (4.0\%) \\  \hline
2 Year & 7  & 17  & 24   & 2   & 2   & 4 \\
&   (1.9\%) &   (1.4\%) &  (1.5\%)  &   (1.1\%) &   (0.3\%) &   (0.5\%)\\  \hline
Other Academic & 24  & 63   &  87   & 16   & 34   & 50 \\ 
&  (6.5\%) &  (5.0\%)  &   (5.4\%) &   (8.7\%) &  (5.0\%) &  (5.8\%)\\ \hline
Research Inst & 10  & 34  & 44 & 4   & 27   & 31   \\ 
&  (2.7\%) &   (2.7\%) &  (2.7\%) &   (2.2\%) &   (4.0\%) &   (3.6\%) \\  \hline
Government & 14   & 48    & 62   & 9  & 7  & 16    \\ 
&   (3.8\%) &  (3.8\%)  &   (3.8\%) &   (4.9\%) &   (1.0\%) &  (1.9\%)  \\ \hline
Business/Industry & 54   & 204  & 258  & 29   & 119   & 148  \\ 
&   (14.6\%) &   (16.3\%) & (15.9\%) &   (15.8\%) &   (17.7\%) &  (17.3\%) \\ \hline
All Others & 96 & 406 & 502& 42 &215 & 257  \\
& (25.4\%) & (31.9\%) & (30.4\%) & (23.0\%) & (32.0\%)& (30.1\%) \\ \hline
Grand Total & 370 & 1250 & 1620 & 183 & 674 & 857 \\
& (100\%)& (100\%)& (100\%)& (100\%)& (100\%)& (100\%) \\ \hline
\end{tabular}
%\end{center}
%\end{minipage}
\end{table}
\medskip

We next consider the most recent cohort of new PhDs from Public Large and Private Large departments and summarize the data in Table \ref{PublicPrivateLargePhDs}.  As usual, the denominators in the percentage of column cohort calculations include new PhDs who accepted other nonacademic jobs,  those who were still seeking or not seeking employment, as well as those whose employment status was unknown at the time of the survey.

We call your attention to a few striking details in Table \ref{PublicPrivateLargePhDs}.  First of all, there is still a very large 
difference in percentages of women and men with PhDs from Public Large and Private Large 
departments who were employed by bachelor's-only departments; the percentage of women taking 
bachelor's-only jobs is at least double the percentage for men. This large gender difference was present in every period of our study and is not diminishing.  Looking at U.S. citizens only, there
 are  considerable  gender differences in PhDs from Private Large programs who were employed by either
  Public or Private Large departments.  For U.S. citizen PhDs from Private Large programs, 
  9.6\% 
  of the women versus 15.5\% of the men were employed by Public Large departments and 12.0\% of the
   women versus 19.1\% of the men took jobs  at Private Large departments.   We will revisit these
    differences in our next section on at least comparable employment rates. We also note that for new
   PhDs from Public Large programs, the percentage of men and women who took jobs in business and 
     industry is higher than the percentage who take jobs at Public Large departments; the same is 
     true for new Private Large PhDs. This was not the case for 1991-2011 PhDs.  Just as in the 
     earlier years, most of the new PhDs who went into business and industry received their 
     degrees from top-producing institutions.  If we look at only U.S. citizens, we see that 19.3\% of
      the women (versus 13.5\% of the men) from Private Large PhD programs took jobs in business 
      and industry.  Looking at U.S. citizen PhDs from Public Large programs, 17.9\% of the men
       and 9.2\% of the women took jobs in business and industry.

\subsection{Statistics and Biostatistics Doctorates:All Departments}

Before we leave this section, we take a brief look and employment patterns for statistics and biostatistics degree recipients. The applied mathematics degree recipients form a much smaller cohort and we will not analyze their employment patterns in this section.   During 1991 and 2000 there were 2,038 degree recipients in statistics or biostatistics (31.9\% women) and between 2001 and 2011 there were 3,582 degree recipients (45.3\% women). Between 2012 and 2015, there were 2,010 degree recipients in statistics or biostatistics (43.4\% women).  It would be interesting to better understand why statistics/biostatistics attracts relatively more women than pure mathematics; 24.8\% of the pure math doctorates between 1991 and 2015 were awarded to women.  The number of statistics/biostatistics degree recipients is increasing over time.  In 2015, we observed that 34\% of all new PhDs wrote dissertations in statistics/biostatistics, more than in any other area. The majority of statistics degree recipients take jobs at departments with doctoral statistics programs, at ``other academic" institutions, in non-US academic departments, in government, or in business and industry.  In the Table \ref{StatPhDs} we will focus our attention to these employment types.  We again recall that in AMS data reports, ``other academic" stands for US academic departments other than pure and applied mathematics departments, statistics and biostatistics departments, departments whose highest degree is Bachelor's or Master's degrees, and 2 year colleges. The All Others row, as well as the denominators in the percentage of column cohort calculations, includes new PhDs who accepted jobs at pure math, applied math, departments whose highest degree is Bachelor's or Master's degrees, and 2 year colleges, those who were still seeking or not seeking employment, and those whose employment status was unknown at the time of the survey.

%\medskip
%
%\begin{minipage}{\textwidth}
%\begin{center}
%Table 3e. Observed Frequencies of First Jobs \\
%(Percentages of Column Totals) for  \\
%Statistics and Biostatistics PhDs 1991--2011
%\medskip
%
%%: Observed Frequencies of First Jobs for Statistics and Biostatisticsonly
%
% \begin{tabular}{|r|c|c|c||c|c|c|} \hline
%& \multicolumn{3}{c||}{1991--2000}&\multicolumn{3}{c|}{2001--2011} \\ \hline \hline
%  Employer Type&Female&Male&Totals&Female&Male&Totals  \\ \hline
%Stat/Biostat & 104 & 202  &306   & 231    & 325 & 556   \\ 
%& (16.0\%) &  (14.6\%) & (15.0\%) & (14.2\%)  & (16.6\%) & (15.5\%) \\ \hline
%Other Academic & 76  & 144   &220   & 294  & 323  & 617   \\ 
%& (11.7\%) &  (10.4\%) &  (10.8\%) &   (18.1\%)  &  (16.5\%) &   (17.2\%) \\ \hline
%
%Research Inst & 37  & 37   & 74  & 108  & 77   & 185   \\ 
%&  (5.7\%) &   (2.7\%) &  (3.6\%) & 
%(6.7\%) &   (3.9\%) &   (5.2\%) \\ \hline
%Government & 35  & 63   & 98  & 99   & 84   & 183  \\ 
%&   (5.4\%) &   (4.5\%) &   (4.8\%) &  (6.1\%) &   (4.3\%) & (5.1\%) \\ \hline
%Business/Industry & 179  &  368  &547 & 408   & 523  & 931   \\ 
%&   (27.5\%) &   (26.5\%) &  (26.8\%) &  (25.1\%) &   (26.7\%) &  (26.0\%)  \\ \hline
%Non US Academic & 57    & 174    & 231   & 
%69  & 104   & 173   \\ 
%&   (8.8\%)  &   (12.5\%)  &  (11.3\%)  & 
%  (4.3\%)&  (5.3\%)  &   (4.8\%) \\ \hline
%
%
%\end{tabular}
%\end{center}
%\end{minipage}
%
%\bigskip

\medskip
%\begin{minipage}{\textwidth}
%\begin{center}
%Table 3e. 
\begin{table}
\centering
\caption{Observed Frequencies of First Jobs 
(Percentages of Column Totals) \\
for Statistics and Biostatistics PhDs 1991--2015}
\label{StatPhDs}
%\end{center}
\medskip

%: Observed Frequencies of First Jobs for Statistics and Biostatistics only
\footnotesize
\begin{tabular*}{1.05\linewidth}{@{\extracolsep{\fill}} |r||c|c|c||c|c|c|c|c|c|}
 %\begin{tabular}[h]{|r||c|c|c||c|c|c|c|c|c|} 
 \hline
& \multicolumn{3}{c||}{1991--2000}&\multicolumn{3}{c|}{2001--2011}& \multicolumn{3}{c|}{2012--2015} \\ \hline \hline
  Employer Type&Female&Male&Totals&Female&Male&Totals &Female&Male&Totals  \\ \hline
Stat/Biostat & 104 & 202  &306   & 231    & 325 & 556  & 101 & 191 & 292 \\ 
& (16.0\%) &  (14.6\%) & (15.0\%) & (14.2\%)  & (16.6\%) & (15.5\%) & (11.6\%) & (16.8\%) & (14.5\%) \\ \hline
Other & 76  & 144   &220   & 294  & 323  & 617 &126 & 141 & 267  \\ 
Academic & (11.7\%) &  (10.4\%) &  (10.8\%) &   (18.1\%)  &  (16.5\%) &   (17.2\%) & (14.4\%) & (12.4\%) & 13.3\%)\\ \hline

Research Inst & 37  & 37   & 74  & 108  & 77   & 185 & 45 & 42 & 87   \\ 
&  (5.7\%) &   (2.7\%) &  (3.6\%) & 
(6.7\%) &   (3.9\%) &   (5.2\%) & (5.2\%) & (3.7\%) & (4.3\%) \\ \hline
Government & 35  & 63   & 98  & 99   & 84   & 183  & 46 & 53 & 99 \\ 
&   (5.4\%) &   (4.5\%) &   (4.8\%) &  (6.1\%) &   (4.3\%) & (5.1\%) & (5.3\%) & (4.7\%) & (4.9\%)\\ \hline
Business/& 179  &  368  &547 & 408   & 523  & 931 & 301 & 411 & 712   \\ 
Industry &   (27.5\%) &   (26.5\%) &  (26.8\%) &  (25.1\%) &   (26.7\%) &  (26.0\%)  & (34.5\%) & (36.1\%) & (35.4\%)\\ \hline
NonUS & 57    & 174    & 231   & 69  & 104   & 173 & 29 & 49 & 78   \\ 
Academic &   (8.8\%)  &   (12.5\%)  &  (11.3\%)  & 
  (4.3\%)&  (5.3\%)  &   (4.8\%) & (3.3\%) & (4.3\%) & (3.9\%) \\ \hline
All Others &180 &510 & 690 & 370 & 505 & 875 & 224 & 251 & 475 \\
&(27.7\%) &(36.7\%) & (33.9\%) & (22.8\%)& (25.8\%)& (24.4\%)& (25.7\%) & (22.1\%) & (23.6\%) \\ \hline
Grand Total & 650& 1388 &2038 & 1623 & 1959 &3582 & 872 & 1138 & 2010 \\
& (100.0\%)& (100.0\%)& (100.0\%)& (100.0\%)& (100.0\%)& (100.0\%)& (100.0\%)& (100.0\%)& (100.0\%) \\ \hline
\end{tabular*}
\end{table}
%\end{minipage}
\normalsize
\medskip

Notice that between 1991 and 2011 a large percentage of statistics/biostatistics degree recipients took jobs in government or in business and industry.  During 2012--2015 the percentage of statistics/biostatistics degree recipients who took jobs in government, business, and industry was even greater.

\section{At Least Comparable Employment Rates}
In this section we ask whether new PhDs attain employment in academia at a level 
at least comparable to that of their degree-granting institution. We consider jobs at 
Research Institutes or Other Non-Profits as desirable and group them with the top-ranking or 
top-producing 
departments (refer to Table~\ref{CompEmpPurePhDs91to11}).\footnote{For Group I PhDs we calculated the percentage who obtained jobs at Group I departments or Research Institutes/Other Non-Profits; for Group II PhDs we calculated the percentage who obtained jobs at Group I - II departments or Research Institutes/Other Non-Profits; and for Group III PhDs we calculated the percentage who obtained jobs at Group I - III departments or Research Institutes/Other Non-Profits.}  As noted in past 
studies \cite{Vitulli Flahive 1997}, \cite{Flahive Vitulli 2010}, since the data collected from departments does not give detailed 
information on the type of position,  a definitive answer to this question is
not possible.  Given that caution,    
 the information in Table~\ref{CompEmpPurePhDs91to11} again indicates 
that women are slightly \textit{less} likely to obtain positions at least comparable with 
their training.  Note that the comparable employment rates for both females and 
males from Group~II institutions  improved between 2001 and 2011,  
particularly the rates for women.  The comparable employment rate for Public Medium departments between 2012 and 2015 was far less favorable to women.  Doctorates from
Group~I  and Public/Private Large  are most likely to accept employment at a department
comparable to their degree-granting department; the comparable employment rate for doctorates from Private Large institutions  is substantially higher than all other comparable employment rates.  The AMS has recently reported percentages of females produced and hired by the various groups of departments in Supplemental Table F.1: Females as a Percentage of New PhDs Produced and Hired by Doctoral-Granting Department Grouping. The only way to get the most recent table is to download the pdf file of the entire Report on New Doctoral Recipients from the AMS website on the survey \cite{YslasVelez2016} and look at the Supplemental Tables at the end of the file; the Supplemental Tables do not appear in the report that is published in the Notices of the AMS.

\medskip
\small
%\begin{minipage}{\textwidth}
%\begin{center}
\begin{table}
\centering
\caption{At Least Comparable Employment Rates  
for New Pure Math PhDs 1991--2015}
\label{CompEmpPurePhDs91to11}
\smallskip
\small
\begin{tabular*}{1.0\linewidth}{@{\extracolsep{\fill}}|c |c| c| c| c| c| c| c| c| c| c| c|} 
\hline
\multicolumn{12}{|c|}{PhD Granting Institution \T\B} \\ 
\hline
\multicolumn{6}{|c|}{1991--2000}& \multicolumn{6}{c|}{  2001--2011 \T\B} \\
\hline
%\hline 

\multicolumn{2}{|c|}{Group I \T\B}& \multicolumn{2}{c|}{Group II \T\B}&
\multicolumn{2}{c|}{Group III \T\B} &\multicolumn{2}{c|}{Group I \T\B}& \multicolumn{2}{c|}{Group II \T\B}&\multicolumn{2}{c|}{Group III \T\B} \\  
\hline
F &  M & F &M & F &  M &F  & M & F&  M & F   & M \\ 
\hline 
22.2\%   &25.3\% & 12.8\% & 13.0\% & 14.3\%  & 11.6\%    &28.5\% & 29.0\% & 18.2\% &  18.9\% & 15.2\%   &17.6\% \\
\hline

\end{tabular*}
\begin{tabular*}{1.0\linewidth}{@{\extracolsep{\fill}}|c|c|c|c|c|c|c|c|c|c|}

\multicolumn{10}{|c|}{2012--2015 \T\B} \\ 

%\multicolumn{10}{|c|}{PhD Granting Institution \T\B} \\ 
\hline

\multicolumn{2}{|c|}{Public Large \T\B }& \multicolumn{2}{c|}{Private Large \T\B}&
\multicolumn{2}{c|}{Public Medium \T\B} &\multicolumn{2}{c|}{Public Small \T\B}& 
\multicolumn{2}{c|}{Private SmalI \T\B} \\  
\hline
F &  M & F &M & F &  M & F  & M & F &  M \\ \hline
21.6\% & 24.7\% & 31.1\% & 32.4\% & 14.7\%&18.6\% &15.8\% &18.8\%&28.1\% &  25.3\%\\
\hline
\end{tabular*}

%\multicolumn{3}{|c|}{Public Large \T\B }& \multicolumn{2}{c|}{Private Large \T\B}&
%\multicolumn{3}{c|}{Public Medium \T\B} &\multicolumn{2}{|c|}{Public Small \T\B}& 
%\multicolumn{2}{c|}{Private Small \T\B} \\  
%\hline
%\multicolumn{2}{|c|}{F} &  M & F &M & \multicolumn{2}{|c|}{F}   &  M & F  & M & F &  M  \\ \hline
%\multicolumn{2}{|c|}{21.6\%} & 24.7\% & 31.1\% & 32.4\% & \multicolumn{2}{|c|}{14.7\%} &18.6\% & 15.8\% &18.8\%&28.1\% &  25.3\% \\
%\hline

%\end{tabular}
\end{table}
%\end{center}
%\end{minipage}
\normalsize
\medskip

During the first period, men from Group I programs were 14.0\% more likely than women to take comparable employment jobs and women from Group III programs were 23.3\% more likely than men to take comparable employment jobs. The gender difference in Group III degree recipients reversed itself in the second period, during which other gender differences were reduced.

Next we turn our attention to the degree recipients between 2012 and 2015 and calculate at least comparable employment rates.
\footnote{For the most recent cohorts of Public and Private Large PhDs we consider employment at Public/Private Large or Research Institutes/Other Non-Profits as employment at least comparable to where the degree was obtained. For Public Medium doctorates, we consider employment at Public/Private Large, Public Medium, or Research Institute/Other Non-Profit as at least comparable employment.  For PhDs from Public or Private Small, we consider employment at Public/Private Large, Public Medium, and Public/Private Small as employment at least comparable to where the degree was obtained.}

\medskip

%\begin{minipage}{\textwidth}
%\begin{center}
%Table 4b. Comparable Employment Rates  \\
%for New Pure Math PhDs 2012--2015\\
%\begin{table}
%\centering
%\caption{Comparable Employment Rates  
%for New Pure Math PhDs 2012--2015}
%\label{CompEmpPurePhDs12to15}
%
%\smallskip
%
%\begin{tabular}[t]{|c|c|c|c|c|c|c|c|c|c|} 
%\hline
%\multicolumn{10}{|c|}{PhD Granting Institution \T\B} \\ 
%\hline
%
%\multicolumn{2}{|c|}{Public Large \T\B }& \multicolumn{2}{|c|}{Private Large \T\B}&
%\multicolumn{2}{|c|}{Public Medium \T\B} &\multicolumn{2}{|c|}{Public Small \T\B}& 
%\multicolumn{2}{|c|}{Private SmalI \T\B} \\  
%\hline
%F &  M & F &M & F &  M & F  & M & F &  M \\ \hline
%21.6\% & 24.7\% & 31.1\% & 32.4\% & 14.7\%&18.6\% &15.8\% &18.8\%&28.1\% &  25.3\%\\
%\hline
%\end{tabular}
%%\end{center}
%%\end{minipage}
%\end{table}
%\medskip

During this period, men from Public Large programs were 14.4\% more likely than women to take
 a job at a department comparable to their degree-granting department. Men from Public Medium
 and Public Small departments were 26.5\%, respectively  19.0\%, more likely than women to take
  a job at a department at least comparable to their degree-granting department.  Women from 
  Private Small departments were 11.1\% more likely than men to take a job at a department at least comparable to their degree-granting department.

 Looking at Table~\ref{CompEmpPurePhDs91to11}
 % and \ref{CompEmpPurePhDs12to15} 
  we see that men who received their degrees from most programs 
were more likely than women to be employed by a department at least comparable to their 
degree-granting department.  The only substantial exception  was the Private Small doctorates during
2012--2015, where women were more likely than men to be employed by a department at least comparable
 to their degree-granting department.

Table~\ref{CompEmpUS91to15} summarizes at least comparable employment rates for U.S. citizens.
 
 \medskip
 \small
% \begin{minipage}{\textwidth}
%\begin{center}
%Table 4c. 
\begin{table}[ht]
\centering
\caption{At Least Comparable Employment Rates for New Pure Math PhDs 
%$\; \; \; \; \; \; \; \; \; \; \; \; \; \; \; \; \; \; \; \; \; \; \; \; \; \; \; \; \; \; \; \; \; \; \; \; \; \; \; \; \; \; \; \; \;$
US Citizens Only   1991--2015}
\label{CompEmpUS91to15}
\smallskip
\small
\begin{tabular*}{1.0\linewidth}{@{\extracolsep{\fill}}|c |c| c| c| c| c| c| c| c| c| c| c|} 
\hline
\multicolumn{12}{|c|}{PhD Granting Institution \T\B} \\ 
\hline
\multicolumn{6}{|c|}{1991--2000}& \multicolumn{6}{c|}{  2001--2011 \T\B} \\
\hline

\multicolumn{2}{|c|}{Group I \T\B}& \multicolumn{2}{c|}{Group II \T\B}&
\multicolumn{2}{c|}{Group III \T\B} &\multicolumn{2}{c|}{Group I \T\B}& \multicolumn{2}{c|}{Group II \T\B}&\multicolumn{2}{c|}{Group III \T\B} \\  
\hline
F &  M & F &M & F &  M &F  & M & F&  M & F   & M \\ 
\hline 
23.4\%   &24.8\% & 12.2\% & 12.4\% & 18.4\%  & 10.9\%    &27.8\% & 30.9\% & 16.5\% &  17.0\% & 16.9\%   &18.8\% \\
\hline
\end{tabular*}
%\end{center}
%\end{minipage}
\begin{tabular*}{1.0\linewidth}{@{\extracolsep{\fill}}|c|c|c|c|c|c|c|c|c|c|}
\multicolumn{10}{|c|}{2012--2105 \T\B} \\ 
\hline

\multicolumn{2}{|c|}{Public Large \T\B }& \multicolumn{2}{c|}{Private Large \T\B}&
\multicolumn{2}{c|}{Public Medium \T\B} &\multicolumn{2}{c|}{Public Small \T\B}& 
\multicolumn{2}{c|}{Private SmalI \T\B} \\  
\hline
F &  M & F &M & F &  M & F  & M & F &  M \\ \hline
21.5\% & 23.7\% & 26.4\% & 38.9\% & 15.6\%&19.4\% &18.3\% &19.3\%&31.6\% &  25.1\%\\
\hline
\end{tabular*}

\end{table}
\normalsize
\medskip

From 1991 to 2000, women from Group III doctoral programs were much more likely to be employed by at least comparable departments than men. The gender difference was substantially greater than when we didn't filter by citizenship.  In both other cases and for all types of  2001--2011 doctorates for U.S. citizens, men were more likely than women to be employed by at least comparable departments.

%\begin{minipage}{\textwidth}
%\begin{center}
%Table 4d. 
%\begin{table}[H]
%\centering
%\caption{Comparable Employment Rates for New Pure Math PhDs 
%US Citizens Only   2012--2015}
%\label{CompEmpUS12to15}
%\smallskip
%
%\begin{tabular}[t]{|c|c|c|c|c|c|c|c|c|c|} 
%\hline
%\multicolumn{10}{|c|}{PhD Granting Institution \T\B} \\ 
%\hline
%
%\multicolumn{2}{|c|}{Public Large \T\B }& \multicolumn{2}{|c|}{Private Large \T\B}&
%\multicolumn{2}{|c|}{Public Medium \T\B} &\multicolumn{2}{|c|}{Public Small \T\B}& 
%\multicolumn{2}{|c|}{Private SmalI \T\B} \\  
%\hline
%F &  M & F &M & F &  M & F  & M & F &  M \\ \hline
%21.5\% & 23.7\% & 26.4\% & 38.9\% & 15.6\%&19.4\% &18.3\% &19.3\%&31.6\% &  25.1\%\\
%\hline
%\end{tabular}
%%\end{center}
%%\end{minipage}
%\end{table}

%\medskip
The striking gender differences for 2012--2015 doctorates were for Private Large and Private Small degree recipients.  For the former group, men were 47.3\% more likely to be employed by at least comparable departments and in the latter group women were 25.9\% more likely to be employed by at least comparable departments. Looking at citizenship differences, women with degrees from Private Large institutions who are U.S. citizens were less likely than non-U.S. citizens to be employed by at least comparable departments; men from Private Large departments who are U.S. citizens were more likely than women to be employed by at least comparable departments.
 
We looked at the question of at least comparable employment in both of our earlier studies. Starting in 
2011, the AMS has reported percentages of females produced and hired by the various groups of 
departments in Supplemental Table F.1: Females as a Percentage of New PhDs Produced and 
Hired by Doctoral-Granting Department Grouping. 
The only way to get the most 
recent table is to download the pdf file of the entire report at the AMS website for the survey \url{http://www.ams.org/profession/data/annual-survey/2015Survey-NewDoctorates-Report.pdf} and look at the Supplemental Tables
at the end of the file.  
 Percentages for males are not reported but can easily be calculated from 
 the table for females.  We hope that in the future AMS will report on at least comparable employment 
 rates.

\section{Conclusion}

We summarize our findings on the four questions that we proposed to study.

\textit{Is the percentage of women (respectively, U.S. citizen) new PhDs increasing?}
The percentage of U.S. citizen women new PhDs increased and then decreased during the three time periods of our study.  Since 2001 non-citizen women received a higher percentage of new degrees than citizen women.  The percentage of non-citizen women receiving PhDs has steadily increased over the three periods of our study.  Non-citizen men received a higher percentage of degrees than citizens during 1991--2000 but the reverse was true in both later periods.

\textit{Are there gender, citizenship or \mbox{gender $\times$ citizenship} differences in initial unemployment rates?} We found that between 1991 and 2000 the unemployment rate 
for non-U.S. citizens was substantially higher than the rate for U.S. citizens. Between 2012 - 2015 the trend reversed itself and the unemployment rate for non-U.S. citizens was considerably lower 
than the rate for citizens. During all periods of our study the unemployment rate for women was 
somewhat lower than the rate for men.  When we looked at \mbox{gender $\times$ citizenship} 
differences we found that between 1991 and 2000 the unemployment rate for women who are 
non-U.S. citizens was more than double the rate for women who are citizens.  

\textit{Are there differences in the type of employment by gender?}  We saw some striking differences in types of employment.  A higher percentage of men take jobs 
 at the top-ranking and top-producing departments than women. A substantially higher percentage of women 
 take jobs at  bachelor's degree 
 only departments  than men.  Moreover,  the percentage of men who 
 take jobs in government, business and industry is considerably higher than the percentage for
  women. 
  
\textit{With regard to academic jobs, 
are men and women equally likely to be employed by departments whose ranking is at least comparable to the degree-granting department?}  During the first period of our study, 1991--2000, the percentage of men receiving degrees from Group I departments who were employed by at least comparable departments was 14.0\% higher than the percentage for women, whereas the percentage of women receiving degrees from Group III departments who were employed by at least comparable departments was 23.3\% higher than the percentage for men.  The gender difference for Group III doctorates reversed itself during 2001-2011, during which time all gender differences were reduced.  Between 2012 and 2015, U.S. citizen males graduating from Private Large departments were much more likely than women to be employed by at least comparable employment departments. During the same time period U.S. citizen females from Private Small departments were more likely than males to be employed by at least comparable departments. Other gender differences in at least comparable employment rates were less pronounced. 

We encourage all doctoral departments and programs to help minimize the number of Unknowns by supplying  as much information about their recent PhDs as possible.  The follow-up survey Employment Experiences of New Doctorates (EENDR) that AMS sends to new PhDs is a less valuable resource since less than half of the new PhDs responded to the survey in recent years \cite{EENDR}. We also encourage new doctorates to return this survey to the AMS after degree completion.

\section{Notes on the Data}
Each year the AMS conducts a census of new PhDs by sending surveys to all departments that 
grant doctoral degrees in mathematics, statistics, applied mathematics and operations research as well as follow-ups to all PhD recipients.  Over the years there 
have been changes in what data is collected and how it is reported.   
Between 1991 and 2011 the AMS reported data for doctorate-granting pure mathematics 
departments partitioned into Groups I, II, and III, based on the latest ranking of U.S. doctoral 
departments as determined by the National Research Council (NRC), a part of the National 
Academies of Science. Starting in 1996 Group I 
was subdivided into Group I Public and Group I Private and Groups IV and V were added.  Group IV consisted of statistics and biostatistics programs and Group V applied mathematics and operations research programs. We  excluded doctorates in operations research (Group Vb) from our current study since they are few in number; during 1991-2011 there were only 174 new PhDs in operations research. In recent years the AMS hasn't surveyed operations research departments or programs. 

The NRC released reports and rankings of research doctoral programs in 1982, 1995
%\footnote{The 
%ranking is based on the 1995 report \textit{Research-Doctorate Programs in the United States: 
%Continuity and Change}, published by The National Academies Press.}, 
and 2010 \cite{NRC1982}, \cite{NRC1995}, \cite{NRC2010}.   Subsequently, the AMS, followed the recommendations of the Joint Data Committee to use these rankings to 
create three groups of pure mathematics doctoral programs, with Group I consisting of the 
top-rated programs.\footnote{Between 1996 and 2011 Group I Public consisted of the top 25 U.S. public 
mathematics departments and Group I Private the top 23 private departments;  Group II contained the next 56 departments;
 Group III contained the remaining U.S. departments reporting a doctoral program in mathematics; 
 Group IV contained U.S. departments (or programs) of statistics, biostatistics, and biometrics 
 reporting a doctoral program; and Group Va consisted of all U.S. departments (or programs) in 
 applied mathematics/applied science reporting a doctoral program.}  The 2010 NRC report \cite{NRC2010} 
% \footnote{\textit{A Data-Based Assessment of Research-Doctorate Programs in the United States}, 
% published by The National Academies Press, \url{https://www.nap.edu/rdp/}}
  does not give a 
 single ranking of programs but rather ranks programs on five different scales with each 
 score presented as a range of rankings; the scales are base on 20 characteristics 
\cite{Mervis 2010}.  Starting in 2012, upon the advice of the Joint Data Committee, 
 the AMS partitioned the pure mathematics departments\footnote{When we speak of  pure mathematics departments we exclude departments in applied mathematics, statistics, and biostatistics} into Math Public Large, Math
 Public Medium, Math Public Small, Math  Private Large, and Math Private Small.  
 This subdivision was 
 based solely on the number of PhDs produced by the departments as reported on the annual 
 surveys between 2000 and 2010.   
% Math 
% Public Large consists of the 25 top-producing public doctoral mathematics programs; Math Public Medium consists of the next 40 public programs and Math Public Small the remaining 64 public 
% departments.  Math Private Large consists of the 24 top-producing private doctoral mathematics
%  programs and Math Private Small the remaining 28 departments. 
  Lists of the departments in 
  these groups as well as a comparison with the former groups can be found on the AMS website for the current annual survey  \cite{AMSSurvey} or past surveys \cite{AMSSurvey_old}.\footnote{As of March 2017, Math Public Large 
  consists of 26 public departments with an annual production rate between 7.0 and 24.2 PhDs per
  year; Math Public Medium consists of 40 public departments with an annual production rate between 3.9 and 6.8 per year; Math Public Small consists of the remaining 64 public departments. Math Private Large consists of 24 private departments with an
   annual production rate between 3.9 and 19.8 PhDs per year and Private Small consists of the remaining 28 departments.}  Due to this change, we do 
   separate analyses for the time periods 1991-2011 and 2012 - 2015.\footnote{ 20 out of 25 former
    Group I Public  and 6 former Group II departments comprise Math Public Large; 22 of the 
    former Group I Private departments plus 1 Group II and 1 Group III departments  comprise 
    Math Private Large; 5 of the former Group I Public departments moved to Math Public Medium;
     they are Georgia Institute of Technology, University of North Carolina at Chapel Hill, University
      of Oregon, University of Utah, and University of Virginia. Arizona State University, Louisiana 
      State University-Baton Rouge, North Carolina State University Texas A\&M University, 
      University of California-Davis, and University of Iowa are the former Group II departments that 
      round out Math Public Large. Claremont Graduate University and Emory University were added to Math Private Large.}
 
The response rate for all groups treated in this report has been very high; the 2015 Annual Survey reports that information was provided by 312 of the 318 doctoral-granting departments queried. 
 Survey response rates by grouping are reported by the AMS in the annual surveys 
 published in the \textit{Notices of the American Mathematical Society} \cite{YslasVelez2016} and available 
 online.  Despite the high overall response rate, over the past several years an increasing number of departments have sent the AMS only basic information on their new PhDs and have often
  omitted data on employment status.  The number of unknowns would be even higher but for web 
  searches by the AMS  that secured additional employment information, especially for those in 
  academia.  This is among the reasons why the AMS conjectures new PhDs who are categorized 
  as Unknowns are skewed toward new PhDs in non-academic employment and individuals who 
  may no longer be in the U.S.  The survey data also do not either distinguish between one-year 
  and multi-year jobs or identify tenure-stream positions.

\section{Acknowledgements}
We thank Thomas H. Barr and Colleen Rose of the American Mathematics 
Society for supplying the data collected from the surveys 
and for answering 
associated questions about the data.

\end{document}